\documentclass[fleqn,11pt,twoside]{article}

\usepackage{amsthm,amsthm,amssymb, color, xcolor,epsfig, graphics, subfigure}

\usepackage{amsmath, graphicx, latexsym, lscape }

\makeatletter
\newcommand{\copyrightnote}[2]{{\renewcommand{\thefootnote}{}
 \footnotetext{\small\it
\begin{flushleft}
 \copyright \ #1   #2  
\end{flushleft}}}}

\newcommand{\Name}[1]{\begin{flushleft}
                       \LARGE \bf #1
                       \end{flushleft}\vspace{-3mm}}

\newcommand{\Author}[1]{\begin{flushleft}
                       \it #1 \end{flushleft}}

\newcommand{\Address}[1]{\begin{flushleft}
                       \it #1 \end{flushleft}}

\newcommand{\Date}[1]{\begin{flushleft}
                      \small  \it #1 \end{flushleft}}

\newcommand{\evenhead}{Author \ name}
\newcommand{\oddhead}{Article \ name}

\renewcommand{\@evenhead}{
\hspace*{-3pt}\raisebox{-15pt}[\headheight][0pt]{\vbox{\hbox to \textwidth
{\thepage \hfil \evenhead}\vskip4pt \hrule}}}
\renewcommand{\@oddhead}{
\hspace*{-3pt}\raisebox{-15pt}[\headheight][0pt]{\vbox{\hbox to \textwidth
{\oddhead \hfil \thepage}\vskip4pt\hrule}}}
\renewcommand{\@evenfoot}{}
\renewcommand{\@oddfoot}{}

\long\def\@makecaption#1#2{%
  \vskip\abovecaptionskip
  \sbox\@tempboxa{\small \textbf{#1.}\ \ #2}%
  \ifdim \wd\@tempboxa >\hsize
    {\small \textbf{#1.}\ \ #2}\par
  \else
    \global \@minipagefalse
    \hb@xt@\hsize{\hfil\box\@tempboxa\hfil}%
  \fi
  \vskip\belowcaptionskip}

\newcommand{\JNMPnumberwithin}[3][\arabic]{%
  \@ifundefined{c@#2}{\@nocounterr{#2}}{%
    \@ifundefined{c@#3}{\@nocnterr{#3}}{%
      \@addtoreset{#2}{#3}%
      \@xp\xdef\csname the#2\endcsname{%
        \@xp\@nx\csname the#3\endcsname .\@nx#1{#2}}}}%
}

\newcommand{\resetfootnoterule} {
  \renewcommand\footnoterule{%
  \kern-3\p@
  \hrule\@width.4\columnwidth
  \kern2.6\p@}
}

\renewcommand{\footnoterule}{}

\makeatother

\theoremstyle{definition}

\begin{document}

\renewcommand{\evenhead}{ {\LARGE\textcolor{blue!10!black!40!green}{{\sf \ \ \ ]ocnmp[}}}\strut\hfill 
C Rogers and A C Briozzo
}
\renewcommand{\oddhead}{ {\LARGE\textcolor{blue!10!black!40!green}{{\sf ]ocnmp[}}}\ \ \ \ \   
Moving boundary problems for a novel extended mKdV equation
}

%%%% Matter for the first page 
\thispagestyle{empty}
\newcommand{\FistPageHead}[3]{
\begin{flushleft}
\raisebox{8mm}[0pt][0pt]
{\footnotesize \sf
\parbox{150mm}{{Open Communications in Nonlinear Mathematical Physics}\ \  \ {\LARGE\textcolor{blue!10!black!40!green}{]ocnmp[}}
\ \ Vol.5 (2025) pp
#2\hfill {\sc #3}}}\vspace{-13mm}
\end{flushleft}}

\FistPageHead{1}{\pageref{firstpage}--\pageref{lastpage}}{ \ \ Article}

\strut\hfill

\strut\hfill

\copyrightnote{The author(s). Distributed under a Creative Commons Attribution 4.0 International License}

\Name{Moving boundary problems for a novel extended mKdV equation. Application of Ermakov-Painlevé II symmetry reduction.}

\Author{Colin Rogers $^{1}$, Adriana C. Briozzo $^{2}$}

\Address{
\small {{$^1$} School of Mathematics and Statistics, The University of New South Wales},\\
\small {Sydney NSW 2052, Australia}\\
\small {{$^2$} Depto. Matem\'atica, FCE, Univ. Austral, Paraguay 1950-CONICET} \\
\small {S2000FZF Rosario, Argentina.}\\
\small{Email: c.rogers@unsw.edu.au,  abriozzo@austral.edu.ar} }

\Date{Received November 6, 2025; Accepted November 18, 2025}

\setcounter{equation}{0}

\smallskip

\noindent
{\bf Citation format for this Article:}\newline
Colin Rogers and Adriana C. Briozzo, 
Moving boundary problems for a novel extended mKdV equation. Application of Ermakov-Painlevé II symmetry reduction,
{\it Open Commun. Nonlinear Math. Phys.}, {\bf 5}, ocnmp:16870, \pageref{firstpage}--\pageref{lastpage}, 2025.

\strut\hfill

\noindent
{\bf The permanent Digital Object Identifier (DOI) for this Article:}\newline
{\it 10.46298/ocnmp.16870}

\strut\hfill

\begin{abstract}
\noindent 
A novel extension of the canonical solitonic mKdV equation is introduced which admits hybrid Ermakov-Painlevé II symmetry reduction. Application of the latter is made to obtain exact solution of Airy-type to a class of moving boundary problems of Stefan kind for this extended mKdV equation. A reciprocal transformation is then applied to the latter to generate an associated exactly solvable class of moving boundary problems for an extension of a base Casimir member of a compacton hierarchy. The extended mKdV equation is shown to be embedded in a range of nonlinear evolution equations with temporal modulation as determined via the action of a class of involutory transformations with origin in Ermakov theory. Associated temporal modulation for the hybrid mKdV and KdV equation as embedded in the classical solitonic Gardner equation is delimited.
\end{abstract}

\label{firstpage}

%%%% The Article text starts here

\section{Introduction}

The mKdV equation has diverse physical applications, notably, in the analysis of nonlinear Alfvén waves in collisionless plasma \cite{kakutani1969} and of acoustic wave propagation in an anharmonic lattice \cite{zabusky1967}. Its connection to the canonical Korteweg-de Vries equation of soliton theory due to Miura \cite{miura1968a} was shown in \cite{rogers1982}  to constitute the spatial part of Bäcklund transformation with classical geometric origins. The property of invariance under Bäcklund transformations and application of iterated associated nonlinear superposition principles has established importance in modern soliton theory \cite{rogers2002}. In a geometric context, the mKdV equation can be derived in connection with the motion of an extensible curve of zero torsion.

Moving boundary problems of Stefan-type in continuum mechanics arise importantly in the analysis of the melting of solids and freezing of liquids (qv. \cite{alexiades1993,crank1987,elliot1982,friedman1982,keller1986,rubinstein1971,tarzia2000} and literature cited therein). The heat balance condition on the moving boundary which separates the phases characteristically provides a nonlinear boundary boundary condition. Reciprocal-type transformations have been applied in \cite{rogers1985} to derive novel analytic solutions to moving boundary problems associated with heat conduction in a range of  metals as detailed by Storm \cite{storm1951} which have temperature-dependent specific heat and thermal conductivity. Conditions for the onset of melting in such metals subject to applied boundary flux may be thereby determined \cite{rogers1988a}. The threshold melting conditions as previously derived by Tarzia \cite{tarzia1981} and Solomon et al \cite{solomon1983} for analogous moving boundary problems for the classical linear heat equation were thereby extended.

In modern soliton theory, a Painlevé II symmetry reduction was applied in \cite{rogers2015} to derive exact solution to a class of moving boundary problems for the canonical Harry Dym equation \cite{vasiliou2001}. The latter arises notably in connection with analysis of the evolution of the interface in Hele-Shaw problems \cite{fokas1998}. A novel sequence of analytically solvable moving boundary problems with interface of the type $x=\gamma t^{n}$ was generated in \cite{rogers2015} via iterated action of a Bäcklund transformation. Exact solutions were derived in terms of Yablonski-Vorob'ov polynomials \cite{vorobav1965,yablonski1959}. Therein the index $n$ adopted a sequence of values of the Painlevé II parameter. In \cite{rogers2017b} such Yablonski-Vorob'ov polynomials solutions for moving boundary problems were shown  to extended to a generalised solitonic Dym equation. The latter was derived in a geometric context in \cite{schief1999} and has physical application to the analysis of peakon solitonic phenomena in hydrodynamics \cite{camassa1993}.

In subsequent work \cite{rogers2022} -\cite{rogers2025b} a series of moving boundary problems of Stefan-type has been shown to be amenable to exact solution via Painlevé II symmetry reduction. In \cite{rogers2025a} it was established that action of a reciprocal transformation links the solitonic Korteweg- de Vries equation to a novel nonlinear evolution equation which incorporates a source term. Moving boundary problems for the latter were thereby shown to admit exact solution. It is remarked that in \cite{briozzo2023} a reciprocal transformation has been applied to reduce to canonical form a class of moving boundary problems involving a source term relevant in a soil mechanics context.

Reciprocal transformations as introduced in a modern solitonic setting in \cite{kingston1982} constitute a class of auto- Bäcklund transformations which act on admitted conservation laws. In \cite{kingston1982} conjugation with the classical geometric nonlinear superposition principle (permutability theorem) of Bianchi was set down which allows iterative generation of multi-soliton solutions in an algorithmic manner.

Reciprocal transformations have been applied in \cite{rogers1984} in the linkage of the canonical AKNS and WKI inverse scattering schemes of \cite{ablowitz1974} and \cite{wadati1979} respectively. Invariance of the 1+1-dimensional Dym solitonic hierarchy under a class of reciprocal transformations was established in \cite{rogers1986a}. Reciprocal transformations in 2+1 dimensions \cite{rogers1986b} were shown in \cite{oevel1993} to connect the Kadomtsev- Petviashvili,  modified Kadomtsev- Petviashvili and 2+1-dimensional Dym triad of S-integrable hierarchies. 

Hybrid Ermaov-Painlevé II systems were originally derived in \cite{rogers2014} via a symmetry reduction of an n+1-dimensional Manakov-type NLS system. Therein, in particular, analysis of certain transverse wave motions in a generalised Mooney- Rivlin hyperelastic material was shown to lead to a novel base canonical Ermakov-Painlevé II reduction. The latter has been subsequently derived and applied in such diverse areas as cold plasma physics \cite{rogers2018a}, Korteweg capillarity theory \cite{rogers2017b} and in connection with Dirichlet-type boundary value problems for a multi-ion Nernst-Planck electrolytic system \cite{amster2015}. A link with the classical Painlevé XXXIV equation was established in \cite{rogers2016b}. Here, a novel extension of the solitonic mKdV equation is introduced which admits Ermakov-Painlevé II symmetry reduction. The latter is applied to derive exact solution to a class of moving boundary problems for this generalised mKdV equation. A reciprocal transformation is then used to derive an extension of a base Casimir member of the compacton hierarchy as set down in \cite{olver1996}. A class of reciprocally associated exactly solvable moving boundary problems for the extended Casimir reciprocal associate is delimited.

\section{An Extended mKdV Equation: Ermakov-Painlevé II Symmetry Reduction} 
Here, a novel mKdV equation with temporal modulation is introduced, namely  
\begin{equation} \label{2.1}
u_{t}-6 u^{2} u_{x} +u_{xxx}+\lambda (t+a)^{\mu}u^{-4}u_x=0 ,\quad \lambda, \mu \in \mathbb{R}
\end{equation} 
and a symmetry reduction to the canonical Ermakov-Painlevé II equation established via the ansatz 
\begin{equation}\label{2.2}
    u=(t+a)^{m} \Psi\bigg(\frac{x}{(t+a)^{n}}\bigg).
\end{equation}
Thus, insertion of the latter into \eqref{2.1}, on reduction, produces
\begin{equation}\label{2.3}
    m\Psi-n\xi\Psi'-6(t+a)^{2m-n+1}\Psi^{2}\Psi'+(t+a)^{-3n+1}\Psi'''+\lambda (t+a)^{\mu-4m-n+1}\Psi^{-4}\Psi'=0,
\end{equation} whence $m=-\frac{1}{3}$, $n=\frac{1}{3}$ together with $\mu=-2$ so that 
\begin{equation}\label{2.4}
   \Psi'''-6\Psi^{2}\Psi'-\frac{1}{3}(\xi\Psi)'+\lambda\Psi^{-4}\Psi'=0
\end{equation}
where, in the preceding, $\xi=x/(t+a)^{n}$.

Integration of the latter yields
\begin{equation}\label{2.5}
   \Psi''-2\Psi^{3}-\frac{1}{3}\xi\Psi-\frac{\lambda}{3}\Psi^{-3}=\zeta, \quad 
\end{equation}
with $\zeta\in\mathbb{R}$. On introduction of the scalings  $\Psi=\delta w^{*}$, $\xi=\epsilon\zeta$ with 
\begin{equation}\label{2.6}
\delta^{2}\epsilon^{2}=1, \quad \epsilon^{3}/3=1
\end{equation}
the canonical Ermakov-Painlevé II equation
\begin{equation}\label{2.8}
    w^{*}_{zz}=2w^{*3}+zw^{*}+\frac{\delta}{w^{*3}}
\end{equation}
results with $\delta=-3\lambda$ if $ \zeta=0$.
\section{A Class of Moving Boundary Problems}
Here a class of Stefan-type moving boundary problems for the extended mKdV equation \eqref{2.1} is introduced, namely 
\begin{equation}\label{3.1}
u_t-6u^{2}u_x+u_{xxx}+\lambda(t+a)^{-2}u^{-4}u_x
=0, \quad 0<x<S(t), \quad t>0,\end{equation}
\begin{equation}\label{3.1.1}
 u_{xx}-2u^3-\frac{\lambda}{3}(t+a)^{-2}u^{-3}=L_m \dot{S}S^i,\;\;\; \text{on}\;\;\ x=S(t),\;\;t>0,
\end{equation}
\begin{equation}\label{3.1.2}
   u=P_mS^j ,\;\;\; \text{on}\;\;\ x=S(t),\;\;t>0,
\end{equation} together with 
\begin{equation}\label{3.2.1}    
u_{xx}-2u^3-\frac{\lambda}{3}(t+a)^{-2}u^{-3}=H_0(t+a)^k,\;\;\; \text{on}\;\;\ x=0,\;\;t>0,
\end{equation}
\begin{equation}\label{3.2.2}
   S(0)=S_0.
   \end{equation}
In the sequel, the moving boundary $x=S(t)=\gamma(t+a)^{1/3}$ is adopted so that the initial condition requires $S_0=\gamma a^{1/3}$.

\subsection*{Boundary conditions}

(I) \begin{equation*}
        u_{xx}(S(t),t)-2u^{3}(S(t),t)-(\lambda/3)(t+a)^{-2}u^{-3}(S(t),t)=L_m S^{i}(t)\dot{S}(t), \quad t>0
    \end{equation*}
Insertion of the symmetry ansatz \eqref{2.2} into the preceding yields
\begin{equation}\label{3.3}
    \Psi''(\gamma)-2\Psi^{3}(\gamma)-(\lambda/3)\Psi^{-3}(\gamma)=L_m S^{i}(t)\dot{S}(t)(t+a)=\frac{1}{3}L_m\gamma^{i+1}(t+a)^{(i+1)/3}
\end{equation}whence $i=-1$ together with 

\begin{equation}\label{3.4}
   \Psi''(\gamma)-2\Psi^{3}(\gamma)-(\lambda/3)\Psi^{-3}(\gamma)=\frac{1}{3}L_m,
\end{equation}so that $L_m=\gamma\Psi(\gamma)$ by virtue of \eqref{2.5} with $\zeta=0$.

(II)
\begin{equation*}
        u(S(t),t)=P_mS^{j}(t), \quad t>0
    \end{equation*}
Accordingly, 
\begin{equation}\label{3.5}
       (t+a)^{-1/3} \Psi(\gamma)=P_m\gamma^{j}(t+a)^{j/3}, \quad t>0
    \end{equation}
so that $j=-1$ and
\begin{equation}\label{3.6}
    P_m=\gamma\Psi(\gamma).
\end{equation}
(III) 

\begin{equation*}
        u_{xx}(0,t)-2u^{3}(0,t)-(\lambda/3)(t+a)^{-2}u^{-3}(0,t)=H_0(t+a)^{k}, \quad t>0.
    \end{equation*}
This yields
\begin{equation}\label{3.7}
       \Psi''(0)-2\Psi^{3}(0)-(\lambda/3)\Psi^{-3}(0)=H_0(t+a)^{k+1}
    \end{equation} so that $k=-1$ and $H_0$ is determined by
\begin{equation}\label{3.8}
       \Psi''(0)-2\Psi^{3}(0)-(\lambda/3)\Psi^{-3}(0)=H_0.
    \end{equation}
Accordingly $H_0=0$ in view of Ermakov-Painlevé II reduction \eqref{2.5} with $\zeta=0$ and $\xi=0$.
\subsection*{\textbf{An Airy reduction}}

The Ermakov-Painlevé II equation \eqref{2.5} with  the alternative scalings $\Psi=\delta w^{*}$, $\xi=\epsilon z$ wherein 
\begin{equation}\label{3.9}
\delta^{2}=-\frac{c_3}{2\epsilon^{2}}, \quad\epsilon^{3}=-3c_2, \quad \sigma=\frac{\lambda}{3}\epsilon^{2}\delta^{-4}
\end{equation} produces symmetry reduction

\begin{equation}\label{3.10}
    w^{*}_{zz}+c_3 w^{*3}+c_2 z w^{*} -\frac{\sigma}{w^{*3}}=\alpha^{*}
\end{equation}
with $\alpha^{*}=\varsigma \epsilon^{2}/\delta\neq0$ if $\zeta\neq0$.

This avatar of the Ermakov-Painlevé II equation has been applied in a Korteweg capillarity system context in \cite{rogers2017b} with zero parameter $\alpha^{*}$. It is remarked that a detailed analysis of aspects of the canonical Painlevé II equation with zero Painlevé parameter was conducted in \cite{bassom1998}.

On setting $w^{*}=\rho^{1/2}$ and with $\alpha^{*}=0$ in \eqref{3.10}, there results
\begin{equation} \label{3.11}
    \rho_{zz}=\frac{1}{2\rho}(\rho_z)^{2}-2c_3 \rho^{2}-2c_2 z\rho+\frac{2\sigma}{\rho}
\end{equation} which with the specifications \cite{rogers2017c}
\begin{equation} \label{3.12}
c_2=-\frac{1}{2}, \quad c_3=-1, \quad \sigma=-\frac{1}{4}\left(\alpha+\frac{1}{2}\right)^{2}
\end{equation} admits the solution
\begin{equation} \label{3.13}
    \rho(z)=w_z+w^{2}+\tfrac{1}{2}z
\end{equation}
whence 
\begin{equation} \label{3.14}
    \Psi=\delta w^{*}=\delta [w_z+w^{2}+\frac{1}{2}z]^{1/2}
\end{equation} with $w(z)$ governed by the classical Painlevé II equation
\begin{equation}\label{3.15}
    w_{zz}=2 w^{3}+zw +\alpha.
\end{equation}
It is remarked that \eqref{3.11} is equivalent to the classical integrable Painlevé XXXIV equation in $\rho$ on appropiate re-scaling.

The latter admits an important particular class of solutions when $\alpha=1/2$, namely
\begin{equation}\label{3.16}
    w=-\frac{\Phi'(z)}{\Phi(z)}
\end{equation} with $\Phi(z)$ governed by de Airy equation
\begin{equation}
    \Phi''+\frac{1}{2}z\Phi=0
\end{equation} whence \eqref{3.14} yields
\begin{equation}
    \Psi=\delta\left[2\left(\frac{\Phi'(z)}{\Phi(z)}\right)^{2}+z/2\right]^{1/2}
\end{equation}
with \begin{equation}
    \Phi=a Ai(-2^{-\frac{1}{3}}z)+b Bi(-2^{-\frac{1}{3}}z),\quad a, b \in \mathbb{R}.
\end{equation}
It is remarked that, in \cite{bass2010} the seed Airy-type solution \eqref{3.16} and the subsequent class of solutions of \eqref{3.15} generated by the iterated action of an admitted Bäcklund transformation has been applied in the solution of certain boundary value problems for the classical Nernst-Planck electrolytic system. Here, use of the Airy-type representation of $\Psi$ in the moving boundary problem determines the parameters $L_m$, $P_m$ and $H_0$ in the boundary conditions.
\section*{An Extended Casimir Equation with Temporal Modulation. A Class of Reciprocal Moving Boundary Problems}
Moving boundary problems of Stefan-type have recently been shown to admit exact solution via Painlevé II symmetry reduction for a range of canonical solitonic equations \cite{rogers2015,rogers2017a,rogers2022,rogers2023,rogers2025a,rogers2025b}. The original investigation of such nonlinear nonlinear moving boundary problems was motivated by aspects of the classical Saffman -Taylor problem with surface tension \cite{saffman1958} and a link to the solitonic Dym equation. Moving boundary problems of a Stefan kind for an extended Dym equation which arises in both hydrodynamics and geometric contexts \cite{schief1999} have subsequently been shown to be solvable by Painlevé II reduction.

Here, a reciprocal transformation is introduced, namely
\begin{equation}\label{3.20}
    dx^{*}=u dx+[-u_{xx}+2u^{3}+\frac{\lambda}{3}(t+a)^{-2}u^{-3}]dt,\quad t^{*}=t, \quad    u^{*}=\frac{1}{u}, \qquad \mathbb{R}^{*}
\end{equation} 
which is applied to the extended mKdV equation \eqref{2.1} and to the class of moving boundary problems \eqref{3.1}-\eqref{3.2.2}. Thus, under $\mathbb{R}^{*}$ there results
\begin{equation}\label{3.21}
    dx=u^{*} dx^{*}+\left[\frac{\partial}{\partial x^{*}}\left(\frac{1}{u^{*}}\frac{\partial}{\partial x^{*}}\right)\left(\frac{1}{u^{*}}\right)-\frac{2}{u^{*2}}-\frac{\lambda}{3}u^{*4}(t^{*}+a)^{-2}\right]dt^{*},
\end{equation} which compatibility condition 

\begin{equation}\label{3.22}
    \frac{\partial u^{*}}{\partial t^{*}}=\frac{\partial}{\partial x^{*}}\left[\frac{\partial}{\partial x^{*}}\left(\frac{1}{u^{*}}\frac{\partial}{\partial x^{*}}\right)\left(\frac{1}{u^{*}}\right)-\frac{2}{u^{*2}}-\frac{\lambda}{3}u^{*4}(t^{*}+a)^{-2}\right],
\end{equation} 
By virtue of the reciprocal connection of the latter via $\mathbb{R}^{*}$ to the extended mKdV equation \eqref{2.1}, it inherits a variant of its Airy-type symmetry reduction.
The reciprocal associate \eqref{3.22} of the extended mKdV equation \eqref{2.1} constitutes a novel extendion of the base Casimir member of the compacton hierarchy of \cite{olver1996}. In \cite{rogers2023} the reciprocal version of the moving boundary problem \eqref{3.1} for the standard Casimir equation corresponding to $\lambda=0$ has been solved via Painlevé II symmetry reduction.
Here, the class of moving boundary problems determined by application of $\mathbb{R}^{*}$ to \eqref{3.1}-\eqref{3.2.2} is determined by the reciprocal system
\[
\tfrac{\partial u^{*}}{\partial t^{*}}=\tfrac{\partial}{\partial x^*}\left[\tfrac{\partial}{\partial x^*}\left(\tfrac{1}{u^*}\tfrac{\partial}{\partial x^*}\left(\tfrac{1}{u^*}\right)\right)-\tfrac{2}{u^{*2}}-\tfrac{\lambda}{3} u^{*4}(t^*+a)^{-2}\right], \quad x^*\big|_{x=0}<x^*<x^*\big|_{x=S(t)}, \quad t^{*}>0,
 \]
 \begin{equation}\label{3.23}
\frac{\partial}{\partial x^*}\left(\frac{1}{u^*}\frac{\partial}{\partial x^*}\left(\frac{1}{u^*}\right)\right)-\frac{2}{u^{*2}}-\frac{\lambda}{3} u^{*4}(t^*+a)^{-2}=L_mS^{i}\dot{S},\quad \text{on} \quad x^{*}|_{x=S(t)},\quad t^{*}>0,
 \end{equation}
 \[\frac{1}{u^{*}}=P_mS^{j}(t),\quad \text{on} \quad x^{*}\big|_{x=S(t)},\quad t^{*}>0,\]
\[\frac{\partial}{\partial x^*}\left(\frac{1}{u^*}\frac{\partial}{\partial x^*}\left(\frac{1}{u^*}\right)\right)-\frac{2}{u^{*2}}-\frac{\lambda}{3} u^{*4}(t^*+a)^{-2}=H_0(t^*+a)^{k},\quad \text{on} \quad x^{*}|_{x=0},\quad t^{*}>0\]
wherein $S(t)=\gamma(t+a)^{1/3}.$ In the preceding 
 \begin{equation}\label{3.24}
 dx^{*}\big|_{x=0}=[-u_{xx}+ 2u^{3}+(\lambda/3)u^{-3}(t+a)^{-2}]dt\big|_{x=0}\end{equation}
 \[=(t+a)^{-1}[-\Psi''+2\Psi^3+(\lambda/3)\Psi^{-3}]dt\big|_{x=0}=0
\]
by virtue of \eqref{2.5} with $\zeta=0$ corresponding to the Ermakov-Painlevé II reduction  \eqref{3.10}. Accordingly, $x^{*}|_{x=0}$ is constant. In addition
 
\[dx^{*}\big|_{x=S(t)}=dx^{*}\big|_{x=\gamma (t+a)^{1/3}}=[(t+a)^{-1/3}\Psi(\gamma)\dot{S}(t)-L_mS^{i}\dot{S}(t)]dt
\]
\begin{equation}\label{3.25}=(t+a)^{-1}[\Psi(\gamma)-L_m\gamma^{-1}]\frac{\gamma}{3}dt
\end{equation}
whence
\begin{equation}\label{3.26}
 x^{*}\big|_{x=S(t)}=\frac{\gamma}{3}[\Psi(\gamma)-L_m\gamma^{-1}]\ln|t^{*}+a|=S^{*}(t^{*})
\end{equation}
upto an additive constant and the reciprocal initial condition on the moving boundary becomes
\begin{equation}
    S^{*}(0)=(\gamma/3)[\Psi(\gamma)-L_m\gamma^{-1}]\ln|a|.
\end{equation}
\section{Modulation}
In \cite{rogers2019}, combined action of a reciprocal and integral-type transformation have been applied sequentially to solve a class of Stefan-type moving boundary problems involving spatial heterogeneity. The latter arise notably as a model of percolation of liquid through porous media in soil mechanics \cite{rogers1988b}. Physical systems which incorporate spatial or temporal modulation arise naturally in both physics and continuum mechanics (qv \cite{belmonte2007,belmonte2008,rogers2016a,zhang2010} and \cite{barclay1978,clements1978,karal1959,rogers2020} respectively together with literature cited therein). Thus, in physics such modulated systems have importance notably in the theory of Bose-Einstein condensates and Bloch wave propagation. In continuum mechanics, modulated systems arise 'inter alia' in elastodynamics \cite{karal1959}, visco-elastodynamics \cite{barclay1978} and the analysis of crack problems in elastostatics \cite{clements1978,rogers2020}.

In recent work \cite{rogers2020}, spatially modulated coupled systems of sine-Gordon, Demoulin and Manakov-type have been systematically reduced to their unmodulated canonical counterparts via classes of involutory transformations. The temporal analogue of the latter to be applied here had their origin in a procedure introduced in \cite{athorne1990}, in connection with the Ermakov-Ray-Reid coupled systems \cite{rogers2012}. The latter has extensive physical applications \cite{rogers2018b}. The transformations are of the type 
\begin{equation}\label{4.1}
    dt^{*}=\rho^{-2}(t)dt, \qquad u^{*}=\rho^{-1}(t) u, \quad\quad T^{*}
\end{equation} and augmented by the relation $\rho^{*}=\rho^{-1}$ admit the key involutory property $T^{**}=I$. Application of $T^{*}$ to \eqref{2.1} results in a wide novel class of extended mKdV equations with temporal modulation, namely
\begin{equation} \label{4.2}
   \frac{\partial}{\partial t^{*}}\left(\frac{u^{*}}{\rho^{*}}\right)-6\rho^{*-5}u^{*2}u_{x*}^{*}+\rho^{*-3}u_{x*x*x*}^{*}+\lambda (t^{*}+a)^{\mu}\rho^{*}u^{*-4}u_{x*}^{*}=0
\end{equation}
with \begin{equation}\label{4.3}
dt=\rho^{*-2}dt^{*} .  
\end{equation}
Application of $T^{*}$ to the moving boundary problems \eqref{3.1}-\eqref{3.2.2} produces associated Stefan-type problems for \eqref{4.2} which inherit the key property of exact solution via Ermakov-Painlevé II symmetry reduction.

In \cite{rogers2020} modulated versions of established solitonic systems were derived via the spatial analogue of the involutory transformation \eqref{4.1}. Therein modulation with $\rho$ determined by hybrid Ermakov-Painlevé II, Ermakov-Painlevé III or Ermakov-Painlevé IV as set down in \cite{rogers2017a} were applied. Here, the modulation term $\rho(t)$ in the class of involutary transformations $T^{*}$ is taken to be determined by the classical Ermakov equation
\begin{equation}\label{4.4}
\rho_{tt}+w(t)\rho=k/\rho^{3}
\end{equation}
which admits the nonlinear superposition principle
\begin{equation}\label{4.5}
    \rho=\sqrt{c_1\Omega_1^{2}+2c_2\Omega_1\Omega_2+c_3\Omega_2^{2}}
\end{equation}
wherein $\Omega_1$, $\Omega_2$ constitute a pair of linearly independent solutions of the auxiliary linear equation
\begin{equation}\label{4.6}
    \Omega_{tt}+w(t)\Omega=0
\end{equation} with constants $k$ together with $c_i$, $i=1,2,3$ such that \begin{equation}\label{4.7}
    c_1c_3-c_2^{2}=\frac{k}{W^{2}}
\end{equation}
where $W=\Omega_1\Omega_{2t}-\Omega_{1t}\Omega_2$ is the constant Wronskian of $\Omega_1,\Omega_2$.
The nonlinear superposition principle \eqref{4.5} can be derived via Lie group invariance as in \cite{rogers1989a}. Therein, application was made to the analysis of moving boundary shoreline evolution hydrodynamics with an underlying rigid basin. In general terms, the classical Ermakov equation \eqref{4.4} has diverse physical applications, notably, 'inter alia', in the nonlinear elastodynamics of boundary-loaded hyperelastic tubes \cite{shahinpoor1971,rogers1989b}, oceanographic pulsrodon eddy evolution \cite{rogers1989c}, magnetogasdynamics \cite{rogers2011a} and the analysis of rotating gas cloud phenomena \cite{rogers2011b}.

\section{An Extended Solitonic Gardner Equation with Temporal Modulation}
A hybrid mKdV and KdV equation, namely the Gardner equation was originally introduced in a now classical paper \cite{miura1968b}. It has subsequently had  diverse physical applications, notably in plasma physics \cite{ruderman2008}, optical lattice theory \cite{wadati1975} and most recently elastodynamics \cite{coclite2021}.

In \cite{slyunyaev1999}, Slyunyaav and Palinovsky established a novel link between the mKdV equation and the canonical Gardner equation
\begin{equation}\label{5.1}
    v_{\tau}+6v(1-v)v_y+v_{yyy}=0.
\end{equation}
Thus, on introduction of 
\begin{equation}\label{5.2}
    u=v-\frac{1}{2}
\end{equation}
together with 
\begin{equation}\label{5.2}
    x=-(3/2)\tau+y, \qquad t=\tau
\end{equation}
the mKdV equation
\begin{equation}\label{5.4}
    u_t-6u^{2}u_x+u_{xxx}=0
\end{equation}
results.

Moving boundary problems of Stefan-type for both the solitonic Gardner equation and a reciprocally associated 3rd order nonlinear evolution incorporating a source term have recently been shown to be amenable to exact solution via a Painlevé II symmetry reduction \cite{rogers2025b}. In the present context, the result of \cite{slyunyaev1999} may be used to establish that the extended Gardner equation
\begin{equation}\label{5.5}
v_{\tau}+6v(1-v)v_y+v_{yyy}+\lambda(\tau+a)^{-2}(v-\frac{1}{2})^{-4}v_y=0
\end{equation}
together with the temporally-modulated class in which it is embedded via application of the involutory transformations $T^{*}$ admit Ermakov-Painlevé II symmetry reduction. Application of the latter can be made to construct exact solution of associated nonlinear moving boundary problems with temporal modulation.

\section{Conclusion}
Investigation of exact solution of Stefan- type moving boundary problems for a range of canonical 1+1- dimensional solitonic equations via Painlevé II symmetry reduction as initiated in \cite{rogers2015} has been detailed in \cite{rogers2017b}, \cite{rogers2022}-\cite{rogers2025b}. In the present work, a class of extended mKdV equations has been introduced which admit solution via hybrid Ermakov-Painlevé II symmetry reduction. Currently under investigation is exact solution of novel
moving boundary problems by means of the latter procedure for certain 2+1- dimensional solitonic equations and their reciprocal associates.
\subsection*{Acknowledgements}

The present work has been partially sponsored by the Project PIP No 11220220100532CO CONICET-UA, Rosario, Argentina

\label{lastpage}
\end{document}